\documentclass[openany]{book}
\usepackage[german]{babel}
\usepackage[all,cmtip]{xy}
\usepackage{hyperref}
\usepackage{amsmath,amssymb,array,fancyhdr,fp,graphicx,
  ifthen,latexsym,multicol,pifont,relsize,stmaryrd,textcomp,yfonts}
\begin{document}
%%%%%%%%%%%%%%%%%%%%%%%%%%%%%%%%%%%%%%%%
\newcommand {\AAM} {\mathbb{A}}
\newcommand {\BB} {\mathbb{B}}
\newcommand {\CC} {\mathbb{C}}
\newcommand {\DD} {\mathbb{D}}
\newcommand {\EE} {\mathbb{E}}
\newcommand {\FF} {\mathbb{F}}
\newcommand {\HH} {\mathbb{H}}
\newcommand {\II} {\mathbb{I}}
\newcommand {\KK} {\mathbb{K}}
\newcommand {\MM} {\mathbb{M}}
\newcommand {\NN} {\mathbb{N}}
\newcommand {\PP} {\mathbb{P}}
\newcommand {\QQ} {\mathbb{Q}}
\newcommand {\RR} {\mathbb{R}}
\newcommand {\TT} {\mathbb{T}}
\newcommand {\YY} {\mathbb{Y}}
\newcommand {\ZZ} {\mathbb{Z}}
\newcommand {\AAA} {\mathcal{A}}
\newcommand {\BBB} {\mathcal{B}}
\newcommand {\CCC} {\mathcal{C}}
\newcommand {\DDD} {\mathcal{D}}
\newcommand {\EEE} {\mathcal{E}}
\newcommand {\FFF} {\mathcal{F}}
\newcommand {\GGG} {\mathcal{G}}
\newcommand {\HHH} {\mathcal{H}}
\newcommand {\III} {\mathcal{I}}
\newcommand {\JJJ} {\mathcal{J}}
\newcommand {\KKK} {\mathcal{K}}
\newcommand {\LLL} {\mathcal{L}}
\newcommand {\MMM} {\mathcal{M}}
\newcommand {\NNN} {\mathcal{N}}
\newcommand {\OOO} {\mathcal{O}}
\newcommand {\PPP} {\mathcal{P}}
\newcommand {\RRR} {\mathcal{R}}
\newcommand {\SSS} {\mathcal{S}}
\newcommand {\TTT} {\mathcal{T}}
\newcommand {\UUU} {\mathcal{U}}
\newcommand {\VVV} {\mathcal{V}}
\newcommand {\XXX} {\mathcal{X}}
\newcommand {\ZZZ} {\mathcal{Z}}
%%%%%%%%%%%%%%%%%%%%%%%%%%%%%%%%%%%%%%%%
\newcommand {\Abo}[1] {\bigskip\NIFG{#1}\m}

\newenvironment {Abstract} [1][Abstract]
{\begin{center}{\Century{1000}\F{#1}}\\[1em]
\begin{minipage}{9cm}\Century{900}} {\end{minipage}\end{center}}

\newcommand {\aitem} [3] {\item [#1] \hypertarget{#2}{} \F{#3:}}

\newcommand {\Align}[1]
{\begin{displaymath}\begin{aligned}#1\end{aligned}\end{displaymath}}

\newenvironment {biblio} {\begin{list}{}
{\renewcommand{\makelabel}[1]{\hfill[##1]}}} {\end{list}}

\newcommand {\F}[1] {\textbf{#1}}
\newcommand {\Frac}[2] {\displaystyle\frac{#1}{#2}}

\newcommand {\gauto} [3][1.5]
{{\renewcommand{\arraystretch}{#1}
\begin{tabular} [t] {|>{$}l<{$}|p{#2cm}|}
\hline&\\[-1.5em] #3 \\[0.5em]\hline\end{tabular}}}

\newcommand {\gk}[1] {\{#1\}}
\newcommand {\hme} {^{-1}}
\newcommand {\K}[1] {{\itshape #1}}

\newcommand {\kapitel} [1] {\setcounter{Kapzaehler}{#1}\setcounter{Satzzaehler}{1}}

\newenvironment {keywords} [1][Keywords]
{\Century{950}\NI\K{Keywords: }} {}

\newcommand {\leer} {=\emptyset}
\newcommand {\M} [1] {{\ensuremath{#1}}}
\newcommand {\m} {\medskip}
\newcommand {\Mal} {\lim\limits}
\newcommand {\Maln}[1][x] {\Mal_{#1 \Mpk 0}}
\newcommand {\Malu}[1][x] {\Mal_{#1 \Mpk \infty}}
\newcommand {\Malun} {\Malu[n]}
\newcommand {\Map} {\prod\limits}
\newcommand {\Mas} {\sum\limits}
\newcommand {\Matrix}[1] {\begin{pmatrix}#1\end{pmatrix}}
\newcommand {\Mf} {\mathop{\bigcirc}\limits}

\newcommand {\MfDreifaelle}[6] {\M{\begin{cases}{#1} & {\quad #2}\\
  {#3} & {\quad #4}\\ {#5} & {\quad #6}\end{cases}}}

\newcommand {\MfZweifaelle}[5][0ex] {\M{\begin{cases}{#2} & {\quad #3}\\[#1]
  {#4} & {\quad #5}\end{cases}}}

\newcommand {\Mk}[3] {\M{\left#1 #3 \right#2}}
\newcommand {\Mkr}[1] {\Mk {(} {)} {#1}}
\newcommand {\mm} {\bigskip}
\newcommand {\Mod} {\operatorname {mod}}
\newcommand {\Movd} {\sqcup}
\newcommand {\Movgd} {\bigsqcup\limits}
\newcommand {\Movg} {\bigcup\limits}
\newcommand {\Mp} {\M{\longrightarrow}}
\newcommand {\Mpk} {\M{\rightarrow}}
\newcommand {\NI} {\noindent}
\newcommand {\NIFG}[1] {\NI\F{\relsize{1}#1}}
\newcommand {\nleer} {\ne\emptyset}

\newcommand {\Nummer} [2][ ]
  {\ifthenelse{\equal{#2}{}}{\Nummernzaehlung[#1].}{\Nummernzaehlung[#1]{} #2.}}

\newcommand {\Nummernzaehlung} [1][ ]
  {\ifthenelse{\value{Kapzaehler}=-1}{}
  {\ifthenelse{\value{Kapzaehler}=0}
  {#1\arabic{Satzzaehler}\stepcounter{Satzzaehler}}
  {#1\arabic{Kapzaehler}.\arabic{Satzzaehler}\stepcounter{Satzzaehler}}}}

\newcommand {\on}[1] {\operatorname{#1}}
\newcommand {\set} {\setlength}
\newcommand {\sm} {\smallskip}
\newcommand {\U}[1] {\underline{#1}}

\newcommand {\zitat} [4][]
{#3\ifthenelse{\equal{#3}{}}{}{ }
[\hyperlink{#2}{\F{#4}}\ifthenelse{\equal{#1}{}}{}{, #1}]}

\newcommand {\zweim} {\hspace*{2em}}
%%%%%%%%%%%%%%%%%%%%%%%%%%%%%%%%%%%%%%%%
 \newcounter{Kapzaehler}
\newcounter{Satzzaehler}
%%%%%%%%%%%%%%%%%%%%%%%%%%%%%%%%%%%%%%%%
\newcommand {\Comment}[2][] {\NI\F{Comment\Nummer{#1}} #2\m}
\newcommand {\Conjecture}[2][] {\NI\F{Conjecture\Nummer{#1}} #2\m}
\newcommand {\Corollary}[2][] {\NI\F{Corollary\Nummer{#1}} \K{#2}\m}
\newcommand {\Definition}[2][] {\NI\F{Definition\Nummer{#1}} #2\m}
\newcommand {\Definitionz}[1][] {\NI\F{Definition\Nummer{#1}}\ }
\newcommand {\Example}[2][] {\NI\F{Example\Nummer{#1}} #2\m}
\newcommand {\Examples}[2][] {\NI\F{Examples\Nummer{#1}} #2\m}
\newcommand {\Examplesz}[1][] {\NI\F{Examples\Nummer{#1}}\ }
\newcommand {\Examplez}[1][] {\NI\F{Example\Nummer{#1}}\ }
\newcommand {\Exercise}[2][] {\NI\F{Exercise\Nummer{#1}} #2\m}
\newcommand {\Exercisez}[1][] {\NI\F{Exercise\Nummer{#1}}\ }
\newcommand {\Lemma}[2][] {\NI\F{Lemma\Nummer{#1}} \K{#2}\m}
\newcommand {\Missing}[1][] {\framebox{MISSING {#1}}}
\newcommand {\Notation}[2][] {\NI\F{Notation\Nummer{#1}} #2\m}
\newcommand {\Note}[2][] {\NI\F{Note\Nummer{#1}} #2\m}
\newcommand {\Notez}[1][] {\NI\F{Note\Nummer{#1}}\ }
\newcommand {\Remark}[2][] {\NI\F{Remark\Nummer{#1}} #2\m}
\newcommand {\Remarkz}[1][] {\NI\F{Remark\Nummer{#1}}\ }
\newcommand {\Proof}[1] {\U{Proof.} #1\m}
\newcommand {\Proposition}[2][] {\NI\F{Proposition\Nummer{#1}} \K{#2}\m}
\newcommand {\Reference}[2][] {\NI\F{Reference\Nummer{#1}} #2\m}
\newcommand {\Standing}[2][] {\NI\F{Standing hypothesis\Nummer{#1}} #2\m}
\newcommand {\Test}[2][] {\NI\F{Test\Nummer{#1}} #2\m}
\newcommand {\Testz}[1][] {\NI\F{Test\Nummer{#1}}\ }
\newcommand {\Theorem}[2][] {\NI\F{Theorem\Nummer{#1}} \K{#2}\m}
%%%%%%%%%%%%%%%%%%%%%%%%%%%%%%%%%%%%%%%%
\newlength{\HPT} \set{\HPT}{0.01pt}

\newlength{\HAB} \set{\HAB}{0.0118pt}

\newcommand {\defont}[2] {\newcommand {#1}[1]{\fontfamily{#2}%
\fontsize{##1\HPT}{##1\HAB}\selectfont}}

\defont {\Century}{pnc}

\Century{1100}
\vspace*{2.5cm}

\Abo {Random functions from coupled dynamical systems}

\begin{center}
Lucilla Baldini and Josef Eschgf\"aller\\[2ex]
{\small Universit\`a degli Studi di Ferrara\\[2ex]
lucilla.baldini\symbol{64}unife.it \& esg\symbol{64}unife.it}
\end{center}\mm

\begin{Abstract}
\NI Let $f:T\Mp T$ be a mapping and $\Omega$ be a subset of $T$ which intersects
every (positive) orbit of $f$. Assume that there are given a second dynamical system
$\lambda:Y\Mp Y$ and a mapping $\alpha:\Omega\Mp Y$. For $t\in T$ let
$\delta(t)$ be the smallest $k$ such that $f^k(t)\in\Omega$ and let
$t_\Omega:=f^{\delta(t)}(t)$ be the first element in the orbit of $t$ which belongs to $\Omega$.
Then we define a mapping $F:T\Mp Y$ by $F(t):=\lambda^{\delta(t)}(t_\Omega)$.
\end{Abstract}


\kapitel{1}
\Abo{1. Weak attractors}

\NI We use the notation $\Mf_xf(x)$, introduced in \zitat {7331}{}{7},
for the mapping $x\longmapsto f(x)$. The symbol $\Mf$ can be obtained in Latex with

\verb/\newcommand {\Fun} {\mathop{\bigcirc}\limits}/.\mm

\Standing {Let $T$ be a non-empty set and $f:T\Mp T$\\
a mapping.}

\Definition {For a subset $A\subset T$ we put\sm

$f^*(A) := \gk{t\in T \mid f^k(t)\in A \text{ for some } k\in\NN} =
\Movg_{k=0}^\infty f^{-k}(A)$}

\Definition {A subset $\Omega\subset T$ is called a \K{weak attractor} (of $f$),
if $f^*(\Omega)=T$, i.e., if for every $t\in T$ there exists $k\in\NN$ such that
$f^k(t)\in\Omega$. In this case for $t\in T$ we put
\Align {\delta(t) &:=\delta(t,\Omega,f) :=\min \gk{k\in\NN \mid f^k(t)\in\Omega}\\[1ex]
t_\Omega &:= f^{\delta(t)}(t)}

\NI $t_\Omega$ is therefore the first element of $\Omega$ we reach from $t$ using $f$.}

\Example {Let $T$ be finite. It is well known that then $T$ can be written as
the disjoint union\sm

$T=f^*(M_1)\Movd\ldots\Movd f^*(M_m)$\sm

\NI where $M_1,\ldots,M_m$ are the minimal orbits of the dynamical system $(T,f)$.
On each $M_j$ the restriction $f_{M_j\Mp M_j}$ is a bijection.

A subset $\Omega\subset T$ is a weak attractor iff $\Omega\cap M_j\nleer$
for every $j=1,\ldots,m$.}

\Remark {Let $\Omega$ be a weak attractor and $t\in T\setminus\Omega$.

Then $\delta(f(t))=\delta(t)-1$. Therefore $\delta$ assumes all elements of
$\gk{0,1,\ldots,\delta(t)}$ as values. Furthermore $\Omega=(\delta=0)$.}

\Proposition {Let $\epsilon:T\Mp\NN$ be a function which, for every $t\in T$,
assumes all elements of $\gk{0,1,\ldots,\epsilon(t)}$ as values. Then:\sm

(1) $\Omega:=(\epsilon=0)\nleer$.

(2) There exists a mapping $g:T\Mp T$ such that $\Omega$ is a weak attractor of $g$
and $\epsilon(t)=\delta(t,\Omega,g)$ for every $t\in T$.}

\Proof {Easy. See Prop. 13.9 in \zitat {26000} {} {2}.}

\kapitel{2}
\Abo{2. Coupled dynamical systems}

\Standing {Let the following data be given:\sm

(1) A set $T$ and a mapping $f:T\Mp T$.

(2) A weak attractor $\Omega$ of $f$.

(3) A set $Y$ and a mapping $\lambda:Y\Mp Y$.

(4) A mapping $\alpha:\Omega\Mp Y$.\sm

\NI We suppose that the sets $T$, $\Omega$ and $Y$ be non-empty.

The triple of mappings $(f,\lambda,\alpha)$ is then called a \K{coupled dynamical system}.}

\Remark {In \zitat {26000} {} {2} a slightly more general concept of \K{automatic generator}
(of pseudorandom functions)
has been defined, mainly in order to include also \K{automatic sequences} (as defined
for example in \zitat {24641} {Allouche/Shallit} {1}).}

\Proposition {For $t\in T$ let\sm

$F(t):= \MfZweifaelle {\alpha(t)} {\text{ if } t\in\Omega}
{\lambda(F(f(t)))} {\text{ if } t\notin\Omega}$\sm

\NI In this way we obtain a well defined mapping $F:T\Mp Y$ which extends $\alpha$
and which we call the pseudorandom function generated by the triple $(f,\lambda,\alpha)$.}

\Remark {In Proposition 2.3 for every $t\in T$ we have\sm

$F(t)=\lambda^{\delta(t)}(\alpha(f^{\delta(t)}(t)))
= \lambda^{\delta(t)}(\alpha(t_\Omega))$}

\Remark {If the mapping $\lambda:Y\Mp Y$ is the identity, then
$F(t)=\alpha(t_\Omega)$ for every $t\in T$.}

\Remark {The most studied classical random sequence generators derive from
a mapping $\lambda: Y\Mp Y$ which for every choice of an initial point $y_0\in Y$ gives rise
to a sequence\sm

$F:=\Mf_n\lambda^n(y_0):\NN\Mp Y$\sm

\NI This can (in a trivial way) be considered as a special case of Proposition 2.3: It suffices to
set $T:=\NN$, $\Omega:=\gk{0}$, and\sm

$f(n):=\MfZweifaelle {n-1} {\text{ if } n>0} {0} {\text{ if } n=0}$\sm

\NI with $\alpha:=\Mf_0y_0:\Omega\Mp Y$.\sm

Then $\delta(n)=n$ for every $n\in\NN$ and from Remark 2.4 we have
$F(n)=\lambda^n(y_0)$ for every $n\in\NN$.}

\Remark {In another trivial way one can obtain every mapping\\$F:T\Mp Y$ with the
construction of Proposition 2.3: We put $\Omega:=T$, $\alpha:=F$ with
$f$ and $\lambda$ (both unused) chosen arbitrarily.}

\Remark {Let $a,b\in T\setminus \Omega$ be such that $f(a)=f(b)$.

Then $F(a)=F(b)$.}

\Proof {This follows from Proposition 2.3.}

\Lemma {Let $a,b\in T$ and $j,k\in\NN$ be such that the following conditions hold:\sm

(1) $j\le \delta(a)$ and $k\le\delta(b)$.

\NI\zweim This means that $f^r(a)\notin\Omega$ for\\
\zweim $0\le r<j$ and $f^r(b)\notin\Omega$ for\\\zweim $0\le r<k$.\sm

(2) $j\le k$.\sm

(3) $f^j(a)=f^k(b)$.\sm

\NI Then $F(b)=\lambda^{k-j}(F(a))$.}

\Proof {Let $c:=f^j(a)=f^kb)$ and $m:=\delta(c)$. Then\sm

$\delta(a)=m+j$ and $\delta(b)=m+k$. By condition (3)\sm

$f^{j+m}(a)=f^{k+m}(b)=:\omega\in\Omega$\sm

\NI thus\sm

$F(a)=\lambda^{m+j}(\alpha(f^{m+j}(a)))=\lambda^{m+j}(\alpha(\omega))$\sm

\NI and similarly
\Align {F(b) &= \lambda^{m+k}(\alpha(f^{m+k}(b)))\\
&=\lambda^{m+k}(\alpha(\omega))=\lambda^{k-j+m+j}(\alpha(\omega))\\
&= \lambda^{k-j}(\lambda^{m+j}(\alpha(\omega)))=\lambda^{k-j}(F(a))}}

\Remark {The mapping $\delta:T\Mp\NN$ itself can be considered as a pseudorandom
mapping and as a special case of Proposition 2.3. For this put\sm

$Y:=\NN$, $\lambda:=\Mf_nn+1$, $\alpha:=\Mf_\omega 0:\Omega\Mp\NN$\sm

\NI Then $F(a)=\lambda^{\delta(a)}(0)=\delta(a)$ for
every $a\in T$ and therefore $F=\delta$.}

\kapitel{3}
\Abo{3. Examples}

\Standing {In this section, for each coupled dynamical\\system $(f:T\Mp T, \lambda:Y\Mp Y,
\alpha:\Omega\Mp\Omega)$ we denote by $F$ the function\\generated by the method
of Proposition 2.3.

The examples have been calculated with Pari/GP ([102] und [103] in\\
\zitat {26000} {Baldini} {2}).

For $a\in\NN$ and $b\in\NN+1$ we denote, as in Pari/GP, by $a\setminus b$ the integer
quotient of $a$ by $b$; e.g. $13\setminus 5=2$. The remainder in the division is denoted
by $a\Mod b$.

$\PP$ is the set of primes.}

\Remark {Let $\alpha$ be constant, say $\alpha(\omega)=y_0$ for every $\omega\in\Omega$.
Then\\$F(t)=\lambda^{\delta(t)}(y_0)$ for every $t\in T$.

As in Proposition 1.6 we don't need to specify $f:T\Mp T$ and $\Omega\subset T$ explicitly;
is suffices that there is given a function $\delta:T\Mp\NN$ which for every $t\in T$
assumes all elements of $\gk{0,1,\ldots,\delta(t)}$ as values.}

% DA 15.3.
\Example {Let $T:=\NN$ and $f:\NN\Mp\NN$ be defined by\sm

$f(n)=\MfDreifaelle {n/3} {\text{ if } n\in 3\NN}{(n+1)\setminus2}
{\text{ if } n\in 3\NN+1} {n-2} {\text{ if } n\in 3\NN+2}$\sm

\NI Since $f(n)<n$ for every $n\ge2$, we can choose the set $\Omega:=\gk{0,1,2,3}$
as weak attractor. Define then $Y:=\gk{0,1}$, $\alpha(n):=n\Mod 2$,
$\lambda(y):=1-y$.

The following table gives the values of $F(n)$ for $n=0,\ldots,319$.}\sm

{\topsep=5pt\partopsep=0pt\relsize{-2}\begin{verbatim}
0 1 0 1 1 0 1 0 0 0 1 1 0 1 1 1 1 0 0 0 1 1 0 0 1 0 0 1 0 0 0 0 1 0 1 1 1 1 0 0
0 1 0 1 1 0 1 1 0 1 1 1 1 0 1 1 0 1 1 0 0 1 1 0 0 1 1 0 0 1 0 0 0 0 1 1 1 0 1 1
0 0 0 1 1 0 0 1 0 0 1 0 0 1 0 0 0 0 1 1 0 0 0 0 1 0 1 1 0 0 1 0 1 1 1 0 0 1 1 0
1 0 0 0 0 1 1 1 0 0 0 1 0 1 1 1 1 0 0 1 1 0 1 1 1 1 0 0 0 1 0 0 1 0 1 1 0 0 1 1
1 0 0 1 1 0 0 1 1 1 0 0 1 1 0 1 1 1 1 0 1 1 0 0 1 1 0 1 1 1 1 0 1 1 0 0 0 1 0 1
1 1 1 0 1 1 0 0 0 1 1 0 0 1 0 0 1 1 0 1 0 0 0 0 1 0 0 1 0 1 1 1 1 0 0 0 1 0 1 1
1 1 0 1 1 0 1 1 0 0 0 1 0 0 1 1 1 0 1 1 0 0 0 1 1 0 0 1 0 0 0 0 1 1 1 0 1 0 0 0
0 1 1 0 0 1 0 0 1 0 0 1 1 0 0 1 1 0 0 1 1 1 0 1 0 0 1 0 0 1 0 0 0 1 1 1 0 0 0 0
\end{verbatim}}

\Remark {We shall often subsume the items in Proposition 2.3
within a table of the following form:\sm

\gauto{2}
{T & \ldots\\f(n) & \ldots\\\Omega & \ldots\\Y & \ldots\\
\alpha(n) & \ldots\\\lambda(y) & \ldots}}

% DA 15.4.
\Example {For $n\in\NN+1$ we denote by $\tau(n)$ the number of (positive) divisors
of $n$. Consider then\sm

\gauto{2.5}
{T & $\NN+2$\\f(n) & $\tau(n)$\\\Omega & $\gk{2}$\\Y & $\gk{0,1}$\\
\alpha(n) & $0$\\\lambda(y) & $1-y$}\sm

\NI Observe that $\tau(2)=2$ and $2\le\tau(n)<n$ for $n\ge3$, so that $\Omega:=\gk{2}$
is indeed a weak attractor of $f$.

We calculate the values of $F(n)$ for $2\le n\le 321$:}\sm

{\topsep=5pt\partopsep=0pt\relsize{-2}\begin{verbatim}
0 1 0 1 1 1 1 0 1 1 0 1 1 1 0 1 0 1 0 1 1 1 0 0 1 1 0 1 0 1 0 1 1 1 1 1 1 1 0 1
0 1 0 0 1 1 0 0 0 1 0 1 0 1 0 1 1 1 1 1 1 0 0 1 0 1 0 1 0 1 1 1 1 0 0 1 0 1 0 0
1 1 1 1 1 1 0 1 1 1 0 1 1 1 1 1 0 0 1 1 0 1 0 0 1 1 1 1 0 1 0 1 0 1 0 0 1 1 1 0
1 1 0 1 1 1 0 1 0 1 1 1 1 0 0 1 0 1 1 1 1 1 0 1 1 0 0 1 1 1 0 0 0 1 1 1 1 1 1 1
0 1 0 0 1 1 1 0 0 0 0 1 0 0 0 1 1 1 1 1 0 1 0 1 0 1 0 0 0 1 0 1 1 0 1 1 1 1 1 1
1 1 1 1 1 0 0 1 1 1 0 1 1 1 1 1 1 1 1 1 0 1 1 1 1 1 1 1 0 0 0 1 1 1 0 1 0 1 1 1
0 0 0 0 0 1 0 1 0 1 1 1 1 0 1 1 0 1 1 0 1 1 1 1 0 1 0 1 1 1 0 0 1 0 1 1 1 0 1 1
0 1 0 0 0 1 1 0 0 1 0 1 1 1 0 0 1 1 1 1 1 1 0 1 1 1 1 1 0 1 1 1 1 1 0 1 0 1 0 1
\end{verbatim}}

% DA 15.9.
\Example {For $n\in\NN+1$ we denote by $\sigma(n)$ the sum of
all (positive) divisors of $n$. For $n\ge2$ then $\sigma(n)\ge1+n$, hence the
mapping $\sigma-1:\NN+2\Mp\NN+2$ is well defined. The prime numbers are exactly
the fixed points of $\sigma-1$:\sm

$\PP=\on{Fix}(\sigma-1)$\sm

\NI It is not known whether all orbits of $\sigma-1$ are finite and whether every
orbit ends in a prime. In \zitat [p. 149] {25337} {Guy} {8}, this conjecture is attributed to Erdös;
the final primes of the first orbits are listed on OEIS as sequence A039654.

If we put $T:=\NN+2$, $f:=\sigma-1$ and assume the conjecture to be true,
$\Omega:=\PP$ becomes a weak attractor of $f$.

\gauto{3.5}
{T & $\NN+2$\\f(n) & $\sigma(n)-1$\\\Omega & $\PP$\\Y & $\gk{0,1,2,3,4}$\\
\alpha(n) & $n\Mod 5$\\\lambda(y) & $(3y+2)\Mod 5$}}\m

{\topsep=5pt\partopsep=0pt\relsize{-2}\begin{verbatim}
2 3 2 0 0 2 0 0 3 1 1 3 1 1 2 2 4 4 0 0 1 3 4 2 0 3 2 4 0 1 3 3 1 3 0 2 4 2 4 1
2 3 1 3 0 2 1 2 1 0 3 3 0 0 0 4 4 4 3 1 2 1 2 1 1 2 3 2 1 1 0 3 1 1 4 2 3 4 1 4
3 3 1 3 0 0 4 4 1 2 3 3 1 0 0 2 4 2 4 1 0 3 4 0 2 2 2 4 0 0 3 3 4 1 4 0 4 1 4 0
1 3 1 2 0 2 0 1 0 1 1 2 4 4 4 2 0 4 1 0 0 3 2 4 2 3 4 4 2 1 0 1 0 0 2 2 4 0 4 0
4 3 1 0 0 2 4 0 4 0 3 3 4 3 2 4 4 4 4 1 1 3 4 3 1 0 1 4 4 1 1 3 1 1 4 2 3 4 4 0
3 4 1 0 0 0 1 4 0 1 4 0 4 1 4 2 0 4 1 0 2 3 1 2 0 2 1 4 0 1 4 3 4 0 4 4 0 4 1 1
4 4 1 0 1 2 4 1 3 1 3 0 1 0 2 2 2 1 3 4 4 3 4 4 4 4 4 4 4 1 3 4 4 2 0 2 4 0 4 1
0 3 1 4 1 1 4 2 0 2 2 3 1 4 4 4 4 1 1 0 2 3 4 2 0 2 0 0 0 1 4 3 3 4 1 2 3 4 0 0
\end{verbatim}}
\Lemma {Let $a,b\in\NN+1$ and $a\ne b$. Then\sm

$\max(|a-b|,(a+b)\setminus 2)<\max(a,b)$}

\Proof {We may assume $a>b$. By hypothesis $b>0$, hence \sm

$|a-b|=a-b<a=\max(a,b)$\sm

\NI Also\sm

$(a+b)\setminus2\le \Frac{a+b}{2}<\Frac{a+a}{2}=a=\max(a,b)$}

\Corollary {Define $f:\NN\times\NN\Mp\NN\times\NN$ by\sm

$f(a,b):=\MfDreifaelle{(a,a)} {\text{ if } b=0} {(b,b)} {\text{ if } a=0}
{(|a-b|,(a+b)\setminus2)} {\text{ otherwise}}$\sm

\NI Then the diagonal $\Omega:=\gk{(n,n) \mid n\in\NN}$ is a weak attractor.

Let $Y$ be another set. A map $\alpha:\Omega\Mp Y$ can be identified with a map\\
$\alpha_0:\NN\Mp Y$.

Given $\alpha_0$, for every map $\lambda:Y\Mp Y$ we obtain a map
$F:\NN\times\NN\Mp Y$.}

% DA, 15.7.
\Example {\gauto{3.5}
{T & $\NN\times\NN$\\f(n) & as in Cor. 3.8\\\Omega & $\gk{(n,n)\mid n\in\NN}$\\
Y & $\gk{0,1,2,3}$\\\alpha(n,n) & $(3n+2)\Mod 4$\\\lambda(y) & $(3y+1)\Mod 4$}}

{\topsep=5pt\partopsep=0pt\relsize{-2}\begin{verbatim}
2 0 1 2 3 0 1 2 3 0 1 2 3 0 1 2 3 0 1 2 3 0 1 2 3 0 1 2 3 0 1 2 3 0 1 2 3 0 1 2
0 1 0 1 0 1 0 2 1 1 0 2 0 1 3 1 0 0 2 0 1 1 0 3 2 1 0 1 1 1 3 3 1 2 0 0 2 0 0 1
1 0 0 1 0 2 3 0 1 3 0 0 0 1 0 3 0 0 1 2 0 1 0 1 3 1 0 0 0 0 2 2 0 1 1 3 1 0 3 0
2 1 1 3 0 1 1 0 0 1 0 2 1 3 1 0 2 3 2 1 3 1 1 0 0 2 2 3 2 1 1 1 3 1 0 1 0 0 1 0
3 0 0 0 2 1 3 0 2 1 0 2 3 1 1 0 0 2 1 1 1 0 1 0 2 0 2 0 1 1 3 1 1 0 0 1 0 0 2 0
0 1 2 1 1 1 0 2 0 1 3 0 2 1 0 1 0 2 0 1 3 1 0 0 0 2 3 0 0 0 1 3 0 3 0 0 0 2 2 1
1 0 3 1 3 0 0 1 1 1 3 1 0 1 3 1 0 2 3 0 1 0 3 0 2 2 1 0 1 1 0 1 3 2 3 1 0 0 0 3
2 2 0 0 0 2 1 3 3 0 1 0 0 2 0 1 0 0 2 0 0 1 3 1 0 0 3 2 1 1 3 1 3 0 1 0 0 0 0 1
\end{verbatim}}

\NI We obtain a map $F$ that can be considered as an infinite pseudorandom matrix, which
is symmetric by construction.

The principal diagonal is the periodic sequence $2103 2103 2103 2103 \ldots$.\\
Indeed we have $F(n,n)=\alpha(n,n)=(3n+2)\Mod 4$ for all $n$.\m

\Remark {Let $m\in\NN+2$ and consider a coupled dynamical system\\
$(f:T\Mp T, \lambda: Y\Mp Y, \alpha:\Omega\Mp Y)$
in which $Y=\gk{0,1,\ldots, m-1}$ and\\$\lambda(y)=(y+1)\Mod m$ for every $y\in Y$.
For $0\le j\le m-1$ let $\Omega_j:=\alpha\hme(j)$.

For $k\in\NN$ and $y\in Y$ then $\lambda^k(y)=(y+k)\Mod m$ and therefore\sm

$F(t)=(j+\delta(t))\Mod m$\sm

\NI for $t\in\Omega_j$.

If in particular $\alpha=0$, then $F(t)=\delta(t) \Mod 2$ for every $t\in T$.}

\Example {We show that the Thue-Morse sequence can be obtained by the method
indicated at the end of Remark 3.10. The sequence can be
defined in the following way (cf. \zitat [p. 69] {24558} {Berstel/Karhumäki} {3}):\sm

$x_0=0, \quad x_{2n}=x_n, \quad x_{2n+1}=1-x_n$\sm

\NI Let $T:=\NN, \Omega:=\gk{0}, Y:=\gk{0,1}, \alpha(0)=0, \lambda(y)=1-y$
and define the mapping $f:\NN\Mp\NN$ by\sm

$f(n):=\MfZweifaelle {n+1}{\text { if } n \text{ is even }}{n\setminus 2}{\text{ otherwise}}$

\NI It is clear that then $\Omega$ is a weak attractor of $f$ and that
\Align {\delta(0) &= 0\\\delta(2n) &=\delta(n)+2 \text { for every } n\ge1\\
\delta(2n+1) &=\delta(n)+1 \text{ for every } n\ge0}\sm

\NI and therefore, by Remark 3.10,
\Align{F(0) &= 0\\
F(2n) &=\delta(2n)\Mod 2 = \delta(n)\Mod 2= F(n)\\
F(2n+1) &= (\delta(n)+1) \Mod 2 = 1-\delta(n)\Mod 2 = 1-F(n)}\sm

\NI This shows that $F(n)=x_n$ for every $n\in\NN$.

The Thue-Morse sequence $F$ can therefore be defined by the table\sm

\gauto{3.5}
{T & $\NN$\\f(n) & $n+1$ for $n$ even\\& $n\setminus2$ for $n$ odd\\
\Omega & $\gk{0}$\\Y & $\gk{0,1}$\\\alpha(0) & $0$\\\lambda(y) & $1-y$}}
% DA, 15.15.
\Example {Define $f:\NN+1\Mp\NN+1$ by\sm

$f(n):=\MfZweifaelle{1}{\text{ if \ }n=1} {n-g(n)} {\text{otherwise}}$\sm

\NI where $g(n)$ is the greatest divisor $\ne n$ of $n$. Therefore, if $n$ is prime,
then $g(n)=1$ and $f(n)=n-1$.

The dynamical system $(\NN+1,f)$ has been studied by Collatz (cf. \zitat [p. 241]
{25401} {Lagarias} {10}).
For $n>1$ one has $f(n)<n$, therefore $\Omega:=\gk{1}$ is a weak attractor
of $f$.\m

\gauto{4.5}
{T & $\NN+1$\\f(n) & $1$ if $n=1$\\& $n-g(n)$ otherwise\\
g(n) & greatest  divisor $\ne n$ of $n$\\\Omega & $\gk{1}$\\
Y & $\gk{0,1}$\\\alpha(n) & $0$\\\lambda(y) & $1-y$}}\sm

{\topsep=5pt\partopsep=0pt\relsize{-2}\begin{verbatim}
0 1 0 0 1 1 0 1 0 0 1 0 1 1 1 0 1 1 0 1 0 0 1 1 0 0 0 0 1 0 1 1 1 0 1 0 1 1 1 0
1 1 0 1 1 0 1 0 0 1 1 1 0 1 0 1 0 0 1 1 0 0 0 0 0 0 1 1 1 0 1 1 0 0 0 0 1 0 1 1
0 0 1 0 0 1 1 0 1 0 1 1 1 0 1 1 0 1 1 0 1 0 1 0 1 1 0 0 1 1 1 0 1 1 0 1 1 0 1 0
0 1 1 1 1 1 0 1 0 1 0 1 0 0 1 0 1 0 1 1 1 0 0 0 0 1 0 1 0 1 0 1 1 0 0 1 0 0 0 0
1 1 0 1 0 0 1 1 0 1 0 0 1 0 0 1 1 0 1 1 0 0 0 0 0 0 0 1 0 0 1 0 1 1 0 0 1 0 1 1
1 0 1 1 0 0 1 1 1 0 1 0 1 1 1 1 1 0 0 0 0 0 1 1 0 0 1 0 1 1 1 0 1 0 0 1 1 0 1 1
0 1 0 0 1 0 1 0 1 0 1 0 0 1 0 0 1 1 1 0 1 1 0 0 1 1 1 1 0 0 1 1 1 0 1 1 0 0 1 0
1 0 1 1 1 1 1 1 0 1 0 0 1 1 0 0 1 1 0 0 0 1 1 0 1 0 1 1 1 1 0 0 1 1 1 1 0 1 0 1
\end{verbatim}}

% DA, 15.16.
\Example {Mimicking Example 3.12, we define\sm

$f(n):=\MfZweifaelle{1}{\text{ if\ }n\le1} {n-h(n)} {\text{otherwise}}$\sm

\NI where, for $n\ge2$, $h(n)$ is the smallest divisor $\ne1$ of $n$. Hence $f(p)=0$
if $p$ is prime. We obtain a dynamical system $(\NN,f)$.

For $n>1$ one has $f(n)<n$, therefore $\Omega:=\gk{0,1}$ is a weak attractor.\m

\gauto{4.5}
{T & $\NN$\\f(n) & $n$ if $n\le1$\\& $n-h(n)$ otherwise\\
h(n) & smallest  divisor $\ne 1$ of $n$\\\Omega & $\gk{0,1}$\\
Y & $\gk{0,1}$\\\alpha(n) & $n$\\\lambda(y) & $1-y$}}\sm

{\topsep=5pt\partopsep=0pt\relsize{-2}\begin{verbatim}
0 1 1 1 0 1 1 1 0 0 1 1 0 1 1 1 0 1 1 1 0 0 1 1 0 1 1 1 0 1 1 1 0 0 1 0 0 1 1 1
0 1 1 1 0 0 1 1 0 0 1 1 0 1 1 0 0 0 1 1 0 1 1 1 0 1 1 1 0 0 1 1 0 1 1 1 0 0 1 1
0 0 1 1 0 1 1 1 0 1 1 1 0 0 1 0 0 1 1 1 0 1 1 1 0 0 1 1 0 1 1 1 0 1 1 0 0 0 1 1
0 0 1 1 0 1 1 1 0 0 1 1 0 0 1 1 0 1 1 1 0 0 1 1 0 1 1 1 0 1 1 1 0 0 1 0 0 1 1 1
0 0 1 1 0 0 1 1 0 1 1 1 0 1 1 0 0 0 1 1 0 1 1 1 0 1 1 1 0 0 1 1 0 1 1 1 0 1 1 1
0 0 1 1 0 1 1 1 0 0 1 1 0 0 1 0 0 0 1 1 0 1 1 1 0 0 1 1 0 1 1 1 0 1 1 0 0 0 1 1
0 1 1 1 0 1 1 0 0 0 1 1 0 0 1 1 0 1 1 1 0 0 1 1 0 1 1 1 0 1 1 1 0 0 1 0 0 1 1 1
0 1 1 1 0 0 1 1 0 1 1 1 0 1 1 0 0 0 1 0 0 0 1 1 0 1 1 1 0 0 1 1 0 1 1 1 0 1 1 1
\end{verbatim}}
\Remark {In the following examples we use the $b$-adic representation of
a natural number $n\in\NN+1$, where $b\in\NN+2$. We use the notation\sm

$n=(a_0,a_1,\ldots,a_k)_b=a_0+a_1b+\ldots+a_kb^k$\sm

\NI where we require $a_k\ne0$.

Let $P:\gk{0,1,\ldots,b-1}\Mp\NN$ be a mapping. We define then
$f:\NN+1\Mp\NN+1$ by $f(n):=\Mas_{j=0}^k P(a_j)$. The dynamic properties of this
type of functions have been studied by numerous authors,
in particular by \zitat {25381} {te Riele} {12} and
\zitat {25423} {Stewart} {13}. They are rather complicated, but in
some cases one can find weak attractors, as we shall see in the following examples.}

\Example {In Remark 3.14 assume $b=10$, $P(a)=a^2$, so that\\
$f(n)=\Mas_{j=0}^ka_j^2$.

Porges, in a paper cited in \zitat [p. 374] {25423} {Stewart} {13}
has shown that every orbit ends up in
the fixed point $1$ or in the cycle $(4\ 16\ 37\ 58\ 89\ 145\ 42\ 20)$.

Therefore $\Omega:=\gk{1,37}$ is a weak attractor.

We choose $Y:=\gk{0,1}$, $\alpha(n)=n\Mod 2$ and $\lambda(y):=1-y$
and obtain a sequence $F$ which begins with}\m

{\topsep=5pt\partopsep=0pt\relsize{-2}\begin{verbatim}
1 0 0 1 1 0 0 0 1 0 1 0 1 1 1 0 0 0 1 0 0 1 0 1 0 0 1 0 1 0 1 0 1 1 0 1 1 0 0 1
1 1 1 1 0 1 0 1 1 1 1 0 0 0 0 1 1 0 1 0 0 0 1 1 1 0 1 1 1 0 0 1 1 0 1 1 0 1 0 0
0 0 0 1 0 1 1 1 1 1 1 1 0 1 1 1 0 1 1 0 1 0 1 1 1 0 0 0 1 1 1 1 0 1 0 1 0 1 1 0
1 0 0 1 1 0 1 0 0 1 0 0 0 1 1 0 0 1 0 1 1 1 1 0 0 1 1 1 0 1 0 1 1 0 0 1 0 0 1 0
1 0 0 1 1 0 0 0 0 0 0 1 0 1 0 0 0 0 1 0 1 0 1 1 0 0 0 1 0 1 1 0 0 0 1 0 1 0 1 0
0 1 0 1 0 0 1 0 1 0 1 0 0 1 1 0 1 0 0 1 0 1 1 0 0 0 0 0 0 0 0 1 0 0 1 0 1 1 0 1
\end{verbatim}}
\Lemma {Let $b\in\NN+2$, $P:\gk{0,1,\ldots,b-1}\Mp\NN$ and $f:\NN+1\Mp\NN+1$
be defined as in Remark 3.14.

Then there exists $n_0\in\NN+1$ such that $f(n)<n$ for every $n\ge n_0$.}

\Proof {We follow the proof of Theorem 1 in \zitat [p. 375] {25423} {Stewart} {13}.

Let $M:=\max(P(0),\ldots,P(b-1))$. Since $\Malu[r]\Frac{b^r}{r+1}=\infty$,
there must exist $r_0\in\NN$ such that $\Frac{b^r}{r+1}>M$ for every $r\ge r_0$.

Set $n:=(a_0,\ldots,a_k)_b\in\NN+1$ (with $a_k\ne0$). Then $b^k\le n$, and\sm

$f(n)=\Mas_{j=0}^k P(a_j)\le (k+1)M<b^k\le n$ \quad for $k\ge r_0$, thus
$f(n)<n$ for $k\ge r_0$, that is, for $n\ge b^{r_0}$.}

\Corollary {Let $f$ be defined as in Remark 3.14. Then:

(1) Every orbit of $f$ is finite.

(2) There exists only a finite number of cycles of $f$ and every orbit ends up in exactly one
of these cycles.}

\Proof {(1) Choose $n_0$ as in Lemma 3.16 and let $n\in\NN+2$.

Each time when $f^j(n)\ge n_0$, after a finite number of steps from $f^j(n)$ one arrives
at a value $<n_0$. But the set $\gk{1,\ldots,n_0-1}$ is finite, therefore there have to
exist repetitions in the set $\gk{f^j(n) \mid j\in\NN}$. This implies that the orbit of $n$
is finite.\sm

(2) Every cycle intersects the set $\gk{1,\ldots,n_0-1}$; but the cycles of $f$ are disjoint.}

\Example {Choose $b=10$ and $P(a)=a^4$ in Remark 3.14.
In \zitat {25526} {Chikawa a.o.} {6}
the authors show that the cycles of $f$ are $(1)$, $(1634)$, $(8208)$, $(9474)$,
$(2178\ 6514)$ and $(13139\ 6725\ 4338\ 4514\ 1138\ 4179\ 9219)$.

Therefore we may choose $\Omega:=\gk{1, 1634, 8208, 9474, 2178, 1138}$.

Let $Y$, $\alpha:\Omega\Mp Y$ and $\lambda:Y\Mp Y$ be as in Example 3.15.}\m

{\topsep=5pt\partopsep=0pt\relsize{-2}\begin{verbatim}
1 0 1 1 1 0 1 1 1 0 1 0 0 1 1 1 1 0 1 0 0 0 1 0 1 0 1 1 0 1 0 1 1 0 1 0 0 0 1 1
1 0 0 1 0 1 0 1 1 1 1 1 1 0 1 0 0 1 1 0 1 0 0 1 0 0 0 0 1 1 1 1 0 0 0 0 0 0 0 1
0 1 0 1 1 0 0 1 0 1 1 0 1 1 1 1 0 0 0 0 1 0 0 1 1 1 1 0 1 1 0 1 1 0 0 1 0 0 0 0
1 0 1 0 0 0 1 1 0 0 1 1 1 0 1 0 0 0 0 1 0 0 0 0 0 0 0 0 1 1 0 0 1 0 1 1 0 1 0 1
1 0 0 0 1 0 1 1 1 1 0 1 0 0 0 1 0 0 1 0 0 1 0 0 1 1 0 1 1 1 0 0 0 1 0 1 1 1 0 0
0 0 1 0 1 0 1 1 0 0 1 0 1 0 0 0 1 1 0 0 0 0 0 0 1 0 1 1 0 1 1 0 1 1 0 1 1 0 0 0
0 0 1 0 0 1 1 1 0 1 0 1 0 0 1 0 0 1 0 0 0 0 1 1 0 0 1 1 1 1 1 1 1 1 0 1 1 1 0 1
1 1 0 1 1 1 1 1 0 0 0 0 0 0 0 1 0 0 1 1 0 1 1 0 1 0 0 0 1 0 1 1 1 0 1 0 0 0 0 1
\end{verbatim}}\m
\Example {For $n=(a_0,\ldots,a_k)_{10}\in\NN+1$ let\\
$f(n):=(a_0+a_1+\ldots+a_k)^2$.
\zitat {25527} {Mohanty/Kumar} {11} show that every orbit of $f$ ends up in one of the cycles
$(1)$, $(81)$ and $(169\ 256)$. We may therefore choose the weak attractor
$\Omega:=\gk{1,81,169}$.

Let $Y$, $\alpha:\Omega\Mp Y$ and $\lambda:Y\Mp Y$ be as in Example 3.15.}\m

{\topsep=5pt\partopsep=0pt\relsize{-2}\begin{verbatim}
1 0 1 1 1 0 1 1 1 0 1 0 0 1 1 1 1 0 1 0 0 0 1 0 1 0 1 1 0 1 0 1 1 0 1 0 0 0 1 1
1 0 0 1 0 1 0 1 1 1 1 1 1 0 1 0 0 1 1 0 1 0 0 1 0 0 0 0 1 1 1 1 0 0 0 0 0 0 0 1
0 1 0 1 1 0 0 1 0 1 1 0 1 1 1 1 0 0 0 0 1 0 0 1 1 1 1 0 1 1 0 1 1 0 0 1 0 0 0 0
1 0 1 0 0 0 1 1 0 0 1 1 1 0 1 0 0 0 0 1 0 0 0 0 0 0 0 0 1 1 0 0 1 0 1 1 0 1 0 1
1 0 0 0 1 0 1 1 1 1 0 1 0 0 0 1 0 0 1 0 0 1 0 0 1 1 0 1 1 1 0 0 0 1 0 1 1 1 0 0
0 0 1 0 1 0 1 1 0 0 1 0 1 0 0 0 1 1 0 0 0 0 0 0 1 0 1 1 0 1 1 0 1 1 0 1 1 0 0 0
0 0 1 0 0 1 1 1 0 1 0 1 0 0 1 0 0 1 0 0 0 0 1 1 0 0 1 1 1 1 1 1 1 1 0 1 1 1 0 1
1 1 0 1 1 1 1 1 0 0 0 0 0 0 0 1 0 0 1 1 0 1 1 0 1 0 0 0 1 0 1 1 1 0 1 0 0 0 0 1
\end{verbatim}}\m

\Example {For $n=(a_0,\ldots,a_k)_{10}\in\NN+1$ let\\
$f(n):= \Map_{j=0}^k(a_j+1)$.
In \zitat [p. 342] {25525} {Wagstaff}{14}, it is
shown that every orbit ends up in one of the cycles
$(18)$ and $(2\ 3\ 4\ 5\ 6\ 7\ 8\ 9\ 10)$. We choose $\Omega:=\gk{18,2}$.

Let $Y$, $\alpha:\Omega\Mp Y$ and $\lambda:Y\Mp Y$ be as in Example 3.15.}\m

{\topsep=5pt\partopsep=0pt\relsize{-2}\begin{verbatim}
1 0 0 1 0 1 0 1 0 1 0 0 0 0 1 1 0 0 0 1 0 1 1 0 1 1 1 0 1 0 0 1 0 0 1 1 0 0 0 1
0 0 0 0 1 0 0 0 1 0 1 1 1 1 0 1 1 0 0 1 1 1 1 0 1 0 0 0 1 0 0 1 0 0 1 0 1 0 0 1
1 0 0 0 0 0 0 0 1 0 0 1 0 1 0 1 0 1 0 1 0 0 0 0 1 1 0 1 0 0 0 1 0 0 1 1 0 0 0 0
1 1 1 1 0 1 1 0 0 0 0 1 0 0 1 0 1 0 0 0 0 1 0 1 0 1 0 1 0 1 1 0 1 0 0 1 0 0 1 1
1 1 0 1 1 0 0 0 1 0 0 1 1 0 0 0 1 0 0 1 0 0 0 1 0 0 0 0 0 0 0 0 0 0 1 1 0 0 0 1
0 1 1 0 1 1 1 0 1 0 1 1 1 1 0 1 1 0 0 1 1 0 0 0 0 0 0 0 1 1 1 0 1 0 0 1 0 0 1 0
1 0 0 0 1 0 1 0 0 1 0 0 0 1 0 0 0 0 0 1 1 0 1 0 0 1 1 1 1 1 1 0 0 1 0 1 1 0 1 0
0 0 0 0 0 1 0 1 0 1 0 1 1 0 0 1 1 0 1 0 0 1 0 0 1 1 0 0 0 0 0 1 0 0 1 0 1 0 0 1
\end{verbatim}}
\Lemma {Let $f:\NN\Mp\NN$ be any function.

Then $\Omega:=\gk{n\in\NN \mid f(n)\ge n}$ is a weak attractor of $f$.}

\Proof {Suppose that there exists $n\in\NN$ such that $f^k(n)\notin\Omega$ for
every $k\in\NN$. By definition of $\Omega$ this implies that
\Align {f(n) &<n\\f^2(n)&<f(n)\\f^3(n)&<f^2(n)\\&\cdots}\sm

\NI In this manner we obtain an infinite and
strictly decreasing sequence of natural numbers\sm

$n>f(n)>f^2(n)>f^3(n)>\ldots$

\NI and this is impossible.}

\Corollary {Let $f:\NN\Mp\NN$ be any function and $U:\NN\Mp\NN$ be a mapping
such that $U(n)\le n$ for every $n\in\NN$.

Then $\Omega{'}:=\gk{n\in\NN \mid f(n)\ge U(n)}$ is a weak attractor of $f$.}

\Proof {$f(n)\ge n$ implies $f(n)\ge U(n)$. With $\Omega$ as in Lemma 3.21 we have
$\Omega\subset\Omega{'}$.}

\Corollary {Let $f:\NN\Mp\NN$ be any function and $V:\NN\Mp\NN$ be another
mapping such that $V(n)\ge n$ for every $n\in\NN$.

Then $\Omega{''}:=\gk{n\in\NN \mid V(f(n))\ge n}$ is a weak attractor of $f$.}

\Proof {$f(n)\ge n$ implies $V(f(n))\ge n$. With $\Omega$ as in Lemma 3.21 we have
$\Omega\subset\Omega{''}$.}

\Example {\gauto{3.6}
{T & $\NN$\\f(n) & $(n^n+1)\Mod (2n+1)$\\\Omega & $\gk{n \mid 3f(n)\ge 2n}$\\
Y & $\gk{0,1}$\\\alpha(n) & $n\Mod 2$\\\lambda(y) & $1-y$}}\m

{\topsep=5pt\partopsep=0pt\relsize{-2}\begin{verbatim}
0 1 1 1 0 0 1 1 0 0 1 1 0 1 1 1 0 1 1 1 0 0 1 1 0 1 1 1 0 0 1 1 0 0 0 1 0 1 0 1
0 0 0 1 0 1 0 1 0 1 1 1 0 0 1 1 0 0 1 1 0 1 0 1 0 0 0 1 0 0 0 1 0 0 1 1 0 1 1 1
0 0 0 1 0 1 1 1 0 0 1 1 0 0 0 1 0 1 1 1 0 1 0 1 0 0 0 1 0 1 0 1 1 0 1 1 0 1 1 1
0 0 0 1 0 0 0 1 0 1 1 1 0 1 1 1 0 1 1 1 0 0 1 0 1 1 1 1 0 1 0 1 0 0 1 1 0 1 1 0
0 1 0 1 0 0 1 1 0 1 0 0 0 0 1 1 0 1 0 1 0 1 0 1 0 1 1 0 0 0 1 1 0 1 1 1 0 1 1 1
0 1 0 1 0 1 0 1 0 0 1 1 0 1 0 1 0 1 1 1 0 0 0 1 0 1 0 1 0 1 1 1 0 0 1 1 0 1 0 1
0 0 0 1 0 0 0 0 0 0 0 1 0 1 1 1 0 0 0 1 0 0 0 1 0 0 0 1 0 1 1 1 0 0 0 1 0 1 1 1
0 0 0 1 0 0 0 1 0 1 0 0 0 0 0 1 0 1 0 1 0 1 0 1 0 1 1 1 0 0 0 1 1 1 0 1 0 1 0 1
\end{verbatim}}

\Example {\gauto{7.2}
{T & $\NN+2$\\f(n) & $\Mkr{13n(n-1)(n-2)+\frac{n(n-1)}{2}\Mod(n+1)}$\\
\Omega & $\gk{n \mid f(n)\ge n}$\\Y & $\gk{0,1}$\\
\alpha(n) & $n\Mod 2$\\\lambda(y) & $1-y$}}\m

{\topsep=5pt\partopsep=0pt\relsize{-2}\begin{verbatim}
0 0 1 0 1 1 0 1 1 0 0 1 0 1 1 0 0 0 0 1 0 0 1 1 0 0 0 1 0 1 1 1 1 0 0 1 0 1 1 1
0 0 0 1 1 1 0 1 0 0 0 0 0 0 0 0 0 0 0 0 0 1 1 1 1 0 1 1 0 0 0 0 1 0 1 0 1 1 0 0
0 0 1 1 0 1 1 1 1 1 0 0 1 1 1 1 1 1 0 0 1 0 1 1 1 0 0 0 0 0 1 0 1 0 0 1 1 1 1 0
0 1 1 0 1 1 1 0 1 1 1 0 1 0 1 0 1 1 0 0 0 0 0 1 1 0 1 0 0 1 0 1 0 0 1 1 1 1 0 0
1 1 0 0 0 0 1 1 0 0 0 0 0 1 1 1 0 0 0 1 0 0 1 0 1 1 0 1 0 1 1 1 0 0 0 1 1 1 0 0
0 0 0 1 0 1 0 1 0 0 0 1 0 0 1 0 1 0 1 0 0 1 0 1 1 0 1 1 1 0 0 0 0 0 1 1 0 1 1 1
1 0 0 0 1 0 1 1 1 0 0 1 1 1 1 1 1 1 0 1 0 1 1 1 1 1 0 1 1 1 1 1 0 1 1 1 1 0 1 0
1 0 1 0 1 1 1 0 1 0 0 1 0 1 0 1 1 1 1 0 0 1 0 0 0 1 1 0 1 0 0 1 1 1 0 1 0 1 1 0
\end{verbatim}}
\Proposition {Let $m=2^k\ge4$ and $a,b\in\NN+1$.

With $T:=\gk{0,1,\ldots,2^k-1}$ let $f:T\Mp T$ be defined by\sm

$f(n):=(an+b)\Mod m$\sm

\NI Then the following conditions are equivalent:\sm

(1) $f$ is a cyclic permutation of $T$.\sm

(2) $a\in4\NN+1$ and $b$ is odd.}

\Proof {This is well known, see e.g. \zitat {2129} {Knuth} {9}.}

\Corollary {Let $k\in\NN$ and $f_k:=\Mf_x(5x+1)\Mod 2^k$.

Then $f_k$ is a cyclic permutation of $\gk{0,1,\ldots,2^k-1}$.}

\Example {Since\sm

$\NN+1 = \Movgd_{k=0}^\infty\gk{2^k,2^k+1,\ldots,2^k+2^k-1}$\sm

\NI we can construct a function $f:\NN+1\Mp\NN+1$ by defining\sm

$f(2^k+x):=2^k+((5x+1)\Mod 2^k)$\sm

\NI for $x\in\gk{0,1,\ldots,2^k-1}$. By Corollary 3.27 the orbits of $f$ are exactly the
intervals $\gk{2^k,\ldots,2^k+2^k-1}$ and on each of these intervals $f$ operates as a cyclic
permutation.

We obtain a weak attractor of $f$, if we choose at least one element from each of these
intervals. In particular\sm

$\Omega:=\gk{2^k \mid k\in\NN}$\sm

\NI is a weak attractor.

\gauto{4}
{T & $\NN+1$\\f(2^k+x) & $2^k+((5x+1)\Mod 2^k)$\\\Omega & $\gk{2^k \mid k\in\NN}$\\
Y & $\gk{0,1,2}$\\\alpha(n) & $n\Mod 3$\\\lambda(y) & $(y+1)\Mod 3$}}

{\topsep=5pt\partopsep=0pt\relsize{-2}\begin{verbatim}
1 2 0 1 1 0 2 2 0 1 0 0 2 2 1 1 1 2 2 2 1 0 0 0 2 0 1 1 0 1 2 2 0 1 1 0 0 2 1 1
0 1 0 2 2 2 0 0 2 0 0 1 2 1 2 2 1 2 2 0 1 0 1 1 1 2 0 1 2 0 0 0 1 0 1 0 0 1 1 2
1 2 1 2 0 0 1 1 2 1 1 1 0 2 2 0 2 0 2 2 1 1 2 2 2 2 2 1 1 0 2 1 0 1 2 0 1 2 0 0
0 0 0 2 2 1 0 2 0 0 1 0 0 2 1 1 0 1 2 2 2 0 2 0 0 1 0 1 2 1 0 2 0 2 2 0 2 1 1 1
0 2 0 0 2 2 0 0 0 0 1 2 2 2 1 2 2 0 0 1 2 0 2 1 2 1 1 0 1 0 1 0 2 1 2 2 1 1 2 2
2 2 0 1 1 2 0 1 2 0 2 0 1 2 2 0 1 0 0 2 1 0 0 2 1 1 1 1 0 0 1 1 1 1 0 0 0 1 0 0
1 2 1 2 0 1 1 2 1 2 2 1 0 2 2 1 1 1 0 2 1 0 2 2 1 0 0 1 0 1 0 2 1 0 2 0 1 2 1 0
1 0 1 2 0 0 2 0 1 0 1 1 1 1 2 1 1 2 0 0 0 0 0 1 1 2 2 2 1 1 1 2 1 2 2 1 0 2 2 2
\end{verbatim}}
\Remark {Let $f:\NN\Mp\NN$ be defined by\sm

$f(n):= \MfDreifaelle {1} {\text{ for } n=0} {f(n\setminus2)} {\text{ for } n \text{ odd}}
{f(n/2)+f(n/2-1)} {\text { for } n \text{ even } > 0}$\sm

\NI In \zitat {23742} {Calkin/Wilf} {4}
it is shown that the sequence of quotients $f(n)/f(n+1)$ contains
every rational number $>0$ exactly once .}

\Remark {Let $f:\NN\Mp\NN$ be defined as in Remark 3.29.

Then $f(n)\le n/2$ for every $n\ge5$.}

\Proof {Induction on $n\ge5$.\sm

Since $f(5)=2$, $f(6)=3$, $f(7)=1$, $f(8)=4$, $f(9)=3$, $f(10)=5$, $f(11)=2$, we have
$f(n)\le n/2$ for $5\le n\le 11$, hence in particular for $n=5$.\sm

\U{$n-1\Mp n$:} (i) Assume first that $n\ge 6$ is even.

Then $f(n)=f(n/2)+f(n/2-1)$.

If $n/2-1\ge5$, by induction we have\sm

$f(n)\le\Frac{n}{4}+\Frac{n}{4}-\Frac{1}{2}=\Frac{n}{2}-\Frac{1}{2}\le\Frac{n}{2}$\sm

\NI Otherwise $n=10,8,6$.\sm

(ii) Suppose $n\ge7$ is odd. Then $f(n)=f((n-1)/2)$. If $\Frac{n-1}{2}\ge5$, then by
induction\sm

$f((n-1)/2)\le \Frac{n-1}{4}\le\Frac{n}{2}$\sm

\NI Otherwise $n=11,9,7$.}

\Example {Let $f$ be defined as in Remark 3.29. From Remark 3.30 it follows that
$\Omega:=\gk{1,2}$ is a weak attractor of $f$.

We can therefore define a coupled dynamical system by the table\sm

\gauto{6}
{T & $\NN$\\f(n) & $1$ for $n=0$\\& $f(n\setminus 2)$ for $n$ odd\\
& $f(n/2)+f(n/2-1)$ for $n$ even $>0$\\\Omega & $\gk{0,1}$\\
Y & $\gk{0,1}$\\\alpha(n) & $n\Mod 2$\\\lambda(y) & $1-y$}}

{\topsep=5pt\partopsep=0pt\relsize{-2}\begin{verbatim}
0 1 0 0 1 1 1 0 0 1 0 1 0 1 0 0 0 0 1 1 1 0 1 1 1 0 1 1 1 0 0 0 0 0 0 0 0 1 1 1
0 1 0 0 1 1 0 1 0 1 1 0 0 1 0 1 1 1 0 0 0 0 0 0 1 0 0 0 1 0 0 0 1 0 0 1 1 1 0 1
1 0 0 1 1 0 0 0 1 1 0 1 1 0 0 1 0 0 1 1 0 1 1 0 0 0 1 1 0 0 1 1 0 1 1 1 0 0 1 0
0 0 1 0 0 0 1 0 1 1 0 0 1 0 1 0 0 1 0 0 0 0 1 0 0 1 0 0 1 0 1 1 0 1 0 1 0 0 1 1
1 1 1 0 1 0 0 1 1 1 1 0 1 0 0 0 0 1 1 1 1 0 0 1 0 1 1 0 0 0 0 1 0 0 0 0 1 1 0 1
0 0 1 1 1 1 0 0 0 0 1 0 1 1 1 1 0 0 1 0 1 1 1 1 1 0 0 1 0 1 0 1 1 0 1 0 0 1 0 0
1 0 0 0 0 1 0 0 1 0 1 0 0 1 1 0 0 1 1 1 0 0 0 0 0 1 0 0 0 1 1 0 0 0 1 1 0 0 1 0
1 0 1 0 1 1 1 0 0 0 1 1 0 0 0 0 1 1 0 0 1 1 1 1 1 0 0 1 0 0 0 1 1 0 1 0 1 1 0 1
\end{verbatim}}

\kapitel{4}
\Abo{4. M\"obius transformations on finite fields}

\Lemma {Let $K=GF(q)$ be a finite field of characteristic $\ne2$ and\\
$\alpha\in K\setminus0$.
Then the equation $x^2=\alpha$ has two distinct roots in $L:=GF(q^2)$.}

\Proof {We can choose an algebraic field extension of $K$ in which the equation has a root $\beta$:
$\beta^2=\alpha$.

We must show that $\beta\in L$, i.e., that $\beta^{q^2-1}=1$, being necessarily $\beta\ne0$.

Since $q^2-1$ is even, we have\sm

$\beta^{q^2-1}=(\beta^2)^{\frac{q^2-1}{2}}=\alpha^{\frac{q^2-1}{2}}
=\alpha^{(q+1)(q-1)/2}=\alpha^{2(q-1)/2}=\alpha^{q-1}=1$\sm

\NI The second root is then $-\beta$. Finally $\beta\ne-\beta$, since
$\on{char}K\ne2$.}

\Corollary {Let $K=GF(q)$ be a finite field of characteristic $\ne2$ and\\$b,c\in K$.
Then the equation $x^2+bx+c=0$ has the solutions $x=\Frac{-b\pm\sqrt{b^2-4c}}{2}$,
where by Lemma 4.1 we can calculate the square roots
$\pm\sqrt{b^2-4c}$ in $GF(q^2)$.}

\Proposition {Let $K=GF(q)$ be a finite field of characteristic $\ne2$ and\\
$A=\Matrix{a & b\\c & d}\in GL(2,K)$.

Let $\alpha,\beta\in GF(q^2)$ be the roots of the characteristic polynomial\\
$x^2-(a+d)x+ad-bc$ of $A$. Suppose that the multiplicative order of $\alpha/\beta$ in
$GF(q^2)$ is $q+1$ (this implies that the characteristic polynomial is irreducible).

Define $f:K\Mp K$ by\sm

$f(t) := \MfZweifaelle {\Frac{at+b}{ct+d}} {\text { if } ct+d\ne0}
{a/c} {\text{ otherwise }}$\sm

\NI Then $f$ is a cyclic permutation of $K$.}

\Proof {\zitat [p. 597]
{25558} {\c{C}e\c{s}melio\u{g}lu/W. Meidl/A. Topuzo\u{g}lu} {5}.}

\Remarkz {We use Pari/GP for verifying the hypotheses of Proposition 4.3 for the
matrix $A:=\Matrix{3 & 2\\5 & 1}$ and $K:=GF(1907)$.\m

The characteristic polynomial $x^2-4x-7$ has the roots $\alpha:=2+\sqrt{11}$ and
$\beta:=2-\sqrt{11}$ which must be calculated in $L:=GF(1907^2)$.

First we find a generator $e$ of the field $L$:}

{\topsep=5pt\partopsep=0pt\relsize{-2}\begin{verbatim}
q=1907; q2=q^2; e=ffgen(q2,'e)
\end{verbatim}}\sm

\NI then we find $\alpha$ and $\beta$ with

{\topsep=5pt\partopsep=0pt\relsize{-2}\begin{verbatim}
r=sqrt(11+0*e); alfa=2+r; beta=2-r;
\end{verbatim}}\sm

\NI Now we can verify that the multiplicative order of $\alpha/\beta$ is equal to\\
$1907+1=1908$ using

{\topsep=5pt\partopsep=0pt\relsize{-2}\begin{verbatim}
t_out(fforder(alfa/beta)) \\ 1908
\end{verbatim}}\m

\Example {We can thus apply Proposition 4.3 in order to obtain a
coupled dynamical system:

\gauto{4}
{T & $GF(1907)$\\f(t) & $\Frac{3t+2}{5t+1}$ \quad if $5t+1\ne0$\\& $3/5$ \quad otherwise\\
\Omega & $\gk{1,100,900}$\\Y & $\gk{0,1}$\\\alpha(t) & $0$\\\lambda(y) & $1-y$}}

{\topsep=5pt\partopsep=0pt\relsize{-2}\begin{verbatim}
1 0 0 0 1 1 1 1 1 1 1 0 0 1 1 1 1 0 0 1 0 1 0 1 1 1 0 0 1 0 0 1 0 1 1 0 0 0 1 0
0 0 1 1 1 1 0 1 1 0 1 0 0 1 1 0 0 1 1 0 1 1 0 1 0 1 1 0 0 0 1 1 0 0 0 0 1 0 0 1
1 0 0 0 0 0 1 1 0 1 0 0 0 0 0 1 1 0 1 1 0 0 0 1 0 0 0 1 0 0 0 1 1 0 0 0 0 0 1 0
1 1 1 0 0 0 1 1 1 0 0 1 0 0 1 0 1 0 1 1 0 1 0 0 1 0 1 0 1 0 1 0 0 1 1 1 1 0 0 1
0 0 0 1 1 1 0 0 0 1 0 1 0 0 1 1 1 1 0 0 1 0 1 1 0 1 1 0 0 1 0 0 1 0 1 1 1 1 1 1
0 1 0 1 0 1 0 0 1 0 1 1 0 1 0 1 1 0 1 0 0 1 1 1 1 1 1 0 0 0 1 0 1 0 1 0 1 1 1 1
0 0 0 0 0 0 0 1 1 1 1 0 0 0 1 0 1 0 1 1 1 1 1 0 1 0 0 1 1 1 1 1 1 0 0 0 0 0 0 0
0 1 1 0 1 0 1 1 1 0 0 1 1 1 1 1 1 1 0 0 0 1 0 1 1 1 1 0 0 1 1 0 1 0 1 1 1 0 1 0
\end{verbatim}}\m

\NI Here, after defining the fields $K$ and $L$ as in Remark 4.4, we obtained the
first terms of the sequence $F$ with

{\topsep=5pt\partopsep=0pt\relsize{-2}\begin{verbatim}
f (t) = {my (u);
t=t+0*e; u=5*t+1; if (u, (3*t+2)/u, 3/5)}

inomega (t) = t_pos(t,[1,100,900])

alfa (t) = 0

lam (y) = 1-y

F (t) = if (inomega(t), alfa(t), lam(F(f(t))))

t_fvo(F,[0..799],40)
\end{verbatim}}

\NI using, as in the other examples, the functions
\verb/t_pos/ and \verb/t_fvo/\\from \K{paritools} available on
\href{http://felix.unife.it/++/paritools}{felix.unife.it/++/paritools}.\mm

\NI\F{Concluding remark.} The first author's thesis \zitat {26000}{}{2}
contains more examples, graphical representations,
Fourier transforms and tests.
\newpage
\Abo{References}

\begin{biblio}
\aitem {1} {24641} {J. Allouche/J. Shallit} Automatic sequences. Cambridge UP 2003.
\aitem {2} {26000} {L. Baldini} Analisi armonica e successioni automatiche generalizzate.
Tesi LM Univ. Ferrara 2015.
\aitem {3} {24558} {J. Berstel/J. Karhum\"aki} Combinatorics on words - a tutorial.\\
Bull. EATCS 79 (2003), 178-228.
\aitem {4} {23742} {N. Calkin/H. Wilf} Recounting the rationals.\\
Am. Math. Monthly 107/4 (2000), 360-363.
\aitem {5} {25558} {A. \c{C}e\c{s}melio\u{g}lu/W. Meidl/A. Topuzo\u{g}lu}
On the cycle structure of
permutation polynomials. Finite Fields Appl. 14 (2008), 593-614.
\aitem {6} {25526} {K. Chikawa/K. Iseki/T. Kusakabe} On a problem by H. Steinhaus.\\
Acta Arithmetica 7 (1962), 251-252.
\aitem {7} {7331} {J. Eschgf\"aller} Almost topological spaces.\\
Ann. Univ. Ferrara 30 (1984), 163-183.
\aitem {8} {25337} {R. Guy} Unsolved problems in number theory. Springer 2004.
\aitem {9} {2129} {D. Knuth} The art of computer programming. Volume 2.
Addison-Wesley.
\aitem {10} {25401} {J. Lagarias} The ultimate challenge - the 3x+1 problem.
AMS 2010.
\aitem {11} {25527} {S. Mohanty/H. Kumar} Powers of sums of digits.\\
Math. Magazine 52/5 (1979), 310-312.
\aitem {12} {25381} {H. te Riele} Iteration of number theoretic functions.\\
Nieuw Arch. Wiskunde 1 (1983), 345-360.
\aitem {13} {25423} {B. Stewart} Sums of functions of digits.\\Can. J. Math. 12 (1960), 374-389.
\aitem {14} {25525} {S. Wagstaff} Iterating the product of shifted digits.\\
Fibonacci Quart. 19 (1981), 340-347.
\end{biblio}
%%%%%%%%%%%%%%%%%%%%%%%%%%%%%%%%%%%%%%%%
\end{document}